\documentclass[pdftex,12pt]{amsart}
\usepackage[dvips]{graphicx,color}
\usepackage[pdftex,colorlinks=true, pdfstartview=FitV, linkcolor=blue, citecolor=blue, urlcolor=blue]{hyperref}
\usepackage{a4wide}
\usepackage{amsthm}
\usepackage{amsopn}
\usepackage{amssymb}
\usepackage{amsmath}
\usepackage{latexsym}
\usepackage{graphics}
\usepackage[utf8]{inputenc}
\usepackage[T1]{fontenc}
\usepackage{amsfonts}
\usepackage{color}
\usepackage[all,cmtip]{xy}

\newtheorem{theorem}{Theorem}[section]
\newtheorem{lemma}[theorem]{Lemma}
\newtheorem{corollary}[theorem]{Corollary}
\newtheorem{proposition}[theorem]{Proposition}
\newtheorem{conjecture}{Conjecture}

\newtheorem{mainthm}{Theorem}

\theoremstyle{definition}
\newtheorem{definition}[theorem]{Definition}
\newtheorem{example}[theorem]{Example}

\theoremstyle{remark}
\newtheorem{remark}[theorem]{Remark}

\usepackage{amsfonts}
\usepackage{color}

\newcommand{\ben}{\begin{enumerate}}
\newcommand{\een}{\end{enumerate}}
\newcommand{\ble}{\begin{lemma}}
\newcommand{\ele}{\end{lemma}}
\newcommand{\bth}{\begin{theorem}}
\renewcommand{\eth}{\end{theorem}}
\newcommand{\bpr}{\begin{proposition}}
\newcommand{\epr}{\end{proposition}}
\newcommand{\bco}{\begin{corollary}}
\newcommand{\eco}{\end{corollary}}
\newcommand{\bcon}{\begin{conjecture}}
\newcommand{\econ}{\end{conjecture}}

\newcommand{\bde}{\begin{definition}}
\newcommand{\ede}{\end{definition}}
\newcommand{\bre}{\begin{remark}}
\newcommand{\ere}{\end{remark}}

\newcommand{\bex}{\begin{example}}
\newcommand{\eex}{\end{example}}
\newcommand{\barr}{\begin{array}}
\newcommand{\earr}{\end{array}}
\newcommand{\btab}{\begin{tabular}}
\newcommand{\etab}{\end{tabular}}
\newcommand{\beq}{\begin{equation}}
\newcommand{\eeq}{\end{equation}}
\newcommand{\bea}{\begin{eqnarray*}}
\newcommand{\eea}{\end{eqnarray*}}
\newcommand{\bce}{\begin{center}}
\newcommand{\ece}{\end{center}}
\newcommand{\bpi}{\begin{picture}}
\newcommand{\epi}{\end{picture}}
\newcommand{\bfi}{\begin{figure} \begin{center}}
\newcommand{\efi}{\end{center} \end{figure}}

\newcommand{\bsl}{\begin{slide}{}}
\newcommand{\esl}{\end{slide}}

\newcommand{\pf}{\noindent{\bf Proof}\hspace{7pt}}

\newcommand{\hso}[1]{\hspace{-1pt}}

\newcommand{\into}{\hookrightarrow}
\newcommand{\sbs}{\subset}
\newcommand{\sbe}{\subseteq}

\def\<{\langle}
\def\>{\rangle}

\newcommand{\bC}{{\bf C}}

\newcommand{\bX}{{\bf X}}

\newcommand{\cA}{{\mathcal A}}

\newcommand{\cC}{{\mathcal C}}

\newcommand{\cF}{{\mathcal F}}
\newcommand{\cG}{{\mathcal G}}

\newcommand{\cN}{{\mathcal N}}
\newcommand{\cO}{{\mathcal O}}
\newcommand{\cP}{{\mathcal P}}

\newcommand{\Aut}{\mathop{\rm Aut}\nolimits}
\newcommand{\Inn}{\mathop{\rm Inn}\nolimits}

\newcommand{\Set}{{\mathop{\rm \bf Set}\nolimits}}
\newcommand{\GRep}{\mathop{\rm \bf GRep}\nolimits}
\newcommand{\Group}{\mathop{\rm \bf Group}\nolimits}
\newcommand{\Amal}{\mathop{\rm \bf Amal}\nolimits}

\newcommand{\ad}{\mathop{\rm ad}\nolimits}

\renewcommand{\bar}{\overline}

\newcommand{\id}{\mathop{\rm id}\nolimits}
\newcommand{\im}{\mathop{\rm im}\nolimits}

\newcommand{\xX}{X}
\newcommand{\sS}{\Sigma}

\newcommand{\normal}{\lhd}

\def\flexbox#1{\mathchoice{\mbox{#1}}{\mbox{#1}}{\mbox{\scriptsize #1}}%
{\mbox{\tiny #1}}}
\newcommand{\SL}{{\mathsf{SL}}}

\newcommand{\PG}{\mathop{\flexbox{\rm PG}}\nolimits}

\newcommand{\FF}{{\mathbb F}}

\newcommand{\NN}{{\mathbb N}}

\newcommand{\typ}{\mathop{\rm typ}}%

\newcommand{\obj}{\mathop{\rm obj}\nolimits}

\newcommand{\dfn}{\em}

\newcommand{\after}{\mathbin{ \circ }}

\newcommand{\Stab}{\mathop{\rm Stab}}

\newcommand{\AI}[1]{\item[\rm{(#1)}]}

\renewcommand{\hat}{\widehat}

\newcommand{\vph}{\varphi}

\newcommand{\mn}{\ \medskip \newline }

\renewcommand{\qed}{\hfill $\square$}

\newcounter{romanlistctr}
{\end{list}}%

 \makeatletter \@addtoreset{equation}{section} \makeatother

 \makeatletter
 \def\section{\@startsection {section}{1}{\z@}{-1.5ex plus -.5ex
 minus -.2ex}{1ex plus .2ex}{\large\bf}}

 \makeatletter
 \def\subsection{\@startsection {subsection}{1}{\z@}{-1.5ex plus -.5ex
 minus -.2ex}{1ex plus .2ex}{\bf}}

\def\<{\langle}
\def\>{\rangle}

\newcommand{\ead}{\email}
\newenvironment{keyword}{Keywords: \keywords}{}
\newcommand{\MSC}[1]{MSC #1:\subjclass}
\newcommand{\sep}{\hspace{1ex}}
\begin{document}
\title{$1$-cohomology of simplicial amalgams of groups}
\author{Rieuwert J. Blok}
\address[Rieuwert J. Blok]{department of mathematics and statistics\\
bowling green state university\\
bowling green, oh 43403\\
u.s.a.}
\ead{blokr@member.ams.org}

\author{Corneliu G. Hoffman}
\address[Corneliu G. Hoffman]{university of birmingham\\
edgbaston, b15 2tt\\
u.k.}
\ead{C.G.Hoffman@bham.ac.uk}


\maketitle

\begin{abstract}
We develop a cohomological method to classify amalgams of groups. We generalize this to simplicial amalgams in any concrete category. We compute the non-commutative 1-cohomology for several examples of amalgams defined over small simplices.
\end{abstract}
\paragraph{\bf Acknowledgements}
Part of the work for this paper was completed during an RiP visit to the Matematisches Forschungsinstitut Oberwolfach in Spring of 2011. The authors would like to express their gratitude for the wonderful research environment.

\section{Introduction}
Recognizing the completion $G$ of an amalgam from the multiplication table of that amalgam can be viewed as playing a Sudoku game on the multiplication table of $G$.
More generally, the aim of the game is to decide what $G$ might look like:
You are given a set of subgroups and their intersections and
you need to decide what the largest group containing such a
structure can be. This approach is very useful for example in the
classification of finite simple groups. More precisely,  induction
and local analysis provides a set of subgroups of the minimal
counterexample to the classification and then amalgam type results such as the Curtis-Tits
and Phan theorems show that the group is known after all.

This leaves open the question of whether just the structure of the
subgroups involved determines the group. Most approaches to this
problem~\cite{BeSh2004,Dun2005,Gr2004,Wie2007} use induction together with a lemma by Goldschmidt~\cite{DelGolSte1985,Gol80} that describes the isomorphism
classes of amalgams of two groups in terms of double coset
enumeration. The results presented here are much more general in that they address the classification of amalgams of any finite rank $\ge 2$ and any number of groups.

 In a recent work \cite{BloHof09a} we used Bass-Serre theory of
graphs of groups to classify all possible amalgams of Curtis-Tits
shape with a given diagram. This note describes a method for
higher rank amalgams.
In general an amalgam can be defined over an arbitrary partially
ordered set. In this paper we shall only consider amalgams
defined over the poset of faces of a simplicial complex.  Goldschmidt's Lemma arises as the case where the simplicial complex is an edge on two vertices.

Our starting point is a connected simplicial complex $\xX=(V,\sS)$ and a fixed
amalgam $\cG_0=\{G_\sigma, \psi^\sigma_\tau\mid\sigma\sbe\tau,
\sigma,\tau\in \sS\}$, where the connecting maps $\psi^\sigma_\tau\colon G_\tau\to
G_\sigma$  are injective group homomorphisms whose image we shall denote $\bar{G}_{\sigma,\tau}$. We call such an amalgam a simplicial amalgam.

Note that if  $\cG_0$ is an amalgam, then it is also a complex of groups in the sense of
Bass~\cite{Bas93}, Serre~\cite{Se1980},  and
Haefliger~\cite{Hae1991}. The difference comes from the fact that in
a complex of groups, for each chain of simplices $\sigma\sbe\rho \sbe\tau$ the diagram
  $$\begin{array}{ccc}
  & & G_\sigma\\
  & \nearrow & \\
  G_\rho & & \uparrow \\
  & \nwarrow & \\
  & & G_\tau\\
  \end{array}$$
is required to commute up to inner automorphism of $G_\sigma$, whereas
in a simplicial amalgam, we insist that this inner automorphism be
the identity. As a consequence, there are more complexes of groups
than simplicial amalgams. Thus, the notion of
isomorphism for complexes of groups is weaker than that of
simplicial amalgams.

Our aim is to classify amalgams of type $\cG_0$, where we define an
amalgam of type $\cG_0$ to be an amalgam whose groups $G_{\sigma}$
and $\bar{G}_{\sigma,\tau}$ are those of $\cG_0$ (For precise definitions see Section~\ref{sec:Amalgams and complexes of groups}). Thus the classification reduces to classifying the collections of connecting maps up to isomorphism of the resulting amalgam. To this end we
first create a collection of automorphism groups
$\cA_0=\{A_\sigma, \alpha^\sigma_\tau \mid \sigma\sbe \tau\in \Sigma\}$, where
 $$A_\sigma=\{g\in \Aut(G_\sigma)\mid g(\bar{G}_{\sigma,\tau})=\bar{G}_{\sigma,\tau}\mbox{ for all }\tau\mbox{ with }\sigma\sbe\tau\in\sS\}$$
and, for each pair $(\sigma,\tau)$ with $\sigma\sbe\tau$, we have
a connecting (``restriction'') map $\alpha^\sigma_\tau \colon A_\sigma\to A_\tau$
given by $\ad(\psi^\sigma_\tau)(f)=(\psi^\sigma_\tau)^{-1}\after f\after \psi^\sigma_\tau$. 
For the rest of the paper we will abusively
denote by $(\psi^\sigma_\tau)^{-1}$ (respectively  $(\varphi^\sigma_\tau)^{-1}$ ), the inverse of the isomorphism $\psi^\sigma_
\tau:G_\tau\to \bar{G}_{\sigma, \tau}$ (respectively $\varphi^\sigma_
\tau:G_\tau\to \bar{G}_{\sigma, \tau}$).
We view $\cA_0$ as a coefficient system on the simplicial complex $\xX$.

In Section~\ref{section:cohomology}, we define a non-commutative first cohomology {\em set} $H^1(\xX,\cA_0)$ on $\xX$ with coefficients in $\cA_0$ and in Section~\ref{section:cohomology and amalgams} we use this to prove the following result.

\begin{mainthm}\label{mainthm1}
If $\xX$ is a non-empty connected simplicial complex, then the
isomorphism classes of amalgams of type $\cG_0$ are parametrized by
$H^1(\xX,\cA_0)$.
\end{mainthm}
\noindent The correspondence between $1$-cocycles and amalgams of type $\cG_0$
is constructive, as is the correspondence between $1$-coboundaries and isomorphisms between such amalgams.

In Section~\ref{section:generalizations} we show that there is a result completely analogous to Theorem~\ref{mainthm1} in the much more general setting of simplicial amalgams in any concrete category.

\begin{mainthm}\label{mainthm2}
Let ${}^0\cG$ be a simplicial amalgam over $X$ in a concrete category $\bC$ and let $^0\cA$ be the associated coefficient system.
Then, the isomorphism classes of amalgams of type $^0\cG$ are parametrized by
$H^1(\xX, ^0\!\cA)$.
\end{mainthm}
For definitions and notation see Section~\ref{section:generalizations}. We finish the paper by illustrating the use of Theorem~\ref{mainthm2} with some examples.

\section{Amalgams and complexes of groups}\label{sec:Amalgams and complexes of groups}
\bde\label{dfn:simplicial complex} We define a {\dfn simplicial
complex} to be a pair $\xX=(V,\sS)$ where $V$ is a set of
vertices and $\sS\sbe \cP(V)$ is a collection of finite subsets of
$V$ with the property that $\{v\}\in \sS$ for every $v\in V$ and
if $\tau\in \sS$, then any subset $\sigma\sbe\tau$ also belongs to
$\sS$. An element $\sigma\in \sS$ is called a simplex of rank
$|\sigma|-1$. The boundary $\partial \tau$ of a simplex $\tau$
consists of all simplices of rank $|\tau|-2$ contained in $\tau$.
\ede

From now on we fix a particular connected simplicial complex $\xX=(V,\sS)$,
with $V=\{1,2,\ldots,n\}$ for some $n\in \NN_{\ge 1}$. Given a
simplex $\tau=\{i_1,i_2,\ldots,i_k\}$ with $i_1<i_2<\cdots<i_k$, we
have $\partial\tau=\{\tau_1,\cdots,\tau_k\}$, where
$\tau_j=\tau-\{i_j\}$ for all $j=1,\ldots,k$. The natural ordering
of $V$ now induces an ordering on $\partial \tau$ in which
$\tau_j<\tau_l$, whenever $j<l$. 
We shall write $\bar{\tau}=\tau_1$.

\bde\label{dfn:amalgam} A {\em simplicial amalgam} over the  complex
$\xX=(V,\sS)$ is a collection
$\cG=\{G_\sigma,\varphi^\sigma_\tau\mid \sigma\sbe \tau,
\sigma,\tau\in \sS\}$, where each $G_\sigma$ is a group and, for each
pair $(\sigma,\tau)$ such that $\sigma\sbe \tau$ we have a
monomorphism $\varphi^\sigma_\tau\colon G_\tau\into G_\sigma$,
called an {\em inclusion map} such that, whenever
$\sigma\sbe\rho\sbe\tau$, we have
$\varphi^\sigma_\rho\after\varphi^\rho_\tau=\varphi^\sigma_\tau$.
For simplicity we shall write $
\bar{G}_{\sigma,\tau}=\varphi^\sigma_\tau(G_\tau)\le G_{\sigma}$.
We shall use the shorthand notation $\cG=\{G_\bullet,\varphi_\bullet^\bullet\}$.

A {\em completion} of $\cG$ is a group $G$ together with a
collection  $\phi=\{\phi_\sigma\mid \sigma\in \sS\}$ of
homomorphisms $\phi_\sigma\colon G_\sigma\to G$, such that whenever
$\sigma\sbe \tau$, we have
$\phi_\sigma\after\varphi^\sigma_\tau=\phi_\tau$. The amalgam
$\cG$ is {\em non-collapsing} if it has a non-trivial completion. A
completion $(\hat{G},\hat{\phi})$ is called {\em universal} if for
any completion $(G,\phi)$ there is a  (necessarily unique)
surjective group homomorphism $\pi\colon \hat{G}\to G$ such that
$\phi=\pi\after\hat{\phi}$. \ede

\bde\label{dfn:CT hom} We define a {\dfn homomorphism} between the amalgams
$\cG^{(1)}=\{G_\bullet^{(1)},  {}^{(1)}\varphi_\bullet^\bullet\}$ and
$\cG^{(2)}=\{G_\bullet^{(2)}, {}^{(2)}\varphi_\bullet^\bullet\}$ to be a map
$\phi=\{ \phi_\sigma\mid \sigma\in\sS\}$ where $\phi_\sigma\colon
G_\sigma^{(1)}\to  G^{(2)}_{\sigma}$ are group homomorphisms  such
that \beq\label{eqn:amhom} \phi_{\sigma}\after
{}^{(1)}\varphi^\sigma_\tau ={}^{(2)}\varphi^\sigma_\tau\after
\phi_{\tau}. \eeq We call $\phi$ an {\dfn isomorphism} of amalgams
if $\phi_\sigma$ is bijective for all $\sigma\in \sS$. \ede

\bde\label{dfn:trivial extension functor}
Adopting the notation from Definition~\ref{dfn:CT hom}, suppose $\tilde{X}=(\tilde{V},\tilde{\Sigma})$ is a simplicial complex such that $V\sbe \tilde{V}$ and $\Sigma\sbe \tilde{\Sigma}$.
Given a simplicial amalgam $(\cG_\bullet,\varphi_\bullet^\bullet)$ 
 over $X$, we define a simplicial amalgam $(\tilde{\cG},\tilde{\varphi}_\bullet^\bullet)$ over $\tilde{X}$ as follows.
 $$
 \begin{array}{lcl}
 \tilde{G}_\sigma = \left\{\begin{array}{ll} G_\sigma &\mbox{ if }\sigma\in \Sigma\\
 \{1\} & \mbox{ else }\end{array}\right.
 & \mbox{ and }
 \tilde{\varphi}_\tau^\sigma  = \left\{\begin{array}{ll} \varphi_\tau^\sigma & \mbox{ if }\sigma,\tau\in \Sigma\\
 \tilde{G}_\tau=\{1\}\into \tilde{G}_\sigma & \mbox{ else. }
\end{array}\right.
\end{array}
 $$
Now if $\phi=\{\phi_\sigma\mid \sigma\in \Sigma\}\colon \cG^{(1)}\to\cG^{(2)}$ is a homomorphism, then we define $\tilde{\phi}=\{\tilde{\phi}_\sigma\mid \sigma\in \tilde{\Sigma}\}\colon \tilde{\cG}^{(1)}\to\tilde{\cG}^{(2)}$ as follows.
 $$
 \tilde{\phi}_\sigma  = \left\{\begin{array}{ll} 
 \phi_\sigma & \mbox{ if }\sigma \in \Sigma\\
 \id\colon \tilde{G}_\sigma^{(1)}=\{1\}\to \{1\}=\tilde{G}_\sigma^{(2)} & \mbox{ else. }
\end{array}\right.
$$
\ede
\ble\label{lem:trivial extension functor}
With the notation of Definition~\ref{dfn:trivial extension functor} the assignment $\cG\mapsto \tilde{\cG}$ and $\phi\mapsto \tilde{\phi}$ is an embedding of the category of simplicial amalgams over $X$ into the category of simplicial amalgams over $\tilde{X}$.
In particular, $\phi\colon\cG^{(1)}\stackrel{\cong}{\longrightarrow} \cG^{(2)}$ if and only if $\tilde{\phi}\colon \tilde{\cG}^{(1)}\stackrel{\cong}{\longrightarrow}\tilde{\cG}^{(2)}$.
Moreover, completions are preserved. That is $(G,\phi)$ is a completion of $(\cG_\bullet, \varphi_\bullet^\bullet)$ if and only if 
 $(G,\tilde{\phi})$ (defined in the obvious way) is a completion of $(\tilde{\cG}_\bullet,\tilde{\varphi}_\bullet^\bullet)$. 
\ele
\pf
This is a completely straightforward excercise.
\qed
\mn 
Lemma~\ref{lem:trivial extension functor}
 allows us to replace $X$ by a simplicial complex of rank at least $2$, if necessary. It also allows us to assume that all subsets of $V$ of cardinality $\le k$ are simplices in $\Sigma$.
For the rest of the paper $\cG_0=\{G_\bullet, \psi_\bullet^\bullet\}$
will be a fixed amalgam over $\xX=(V,\sS)$. 
\bde Let $\cG_0=\{G_\bullet, \psi_\bullet^\bullet\}$ be an amalgam over $\xX=(V,\sS)$. An amalgam {\em of type $\cG_0$} is an amalgam $\cG=\{G_\bullet, \varphi_\bullet^\bullet\}$, where, for each $\sigma\in \sS$, the group $G_\sigma$ is that of $\cG_0$ and in which $\varphi^\sigma_\tau(G_\sigma)=\bar{G}_{\sigma,\tau}=\psi^\sigma_\tau(G_\sigma)$ for each pair $(\sigma,\tau)$ with $\sigma\sbe \tau\in \sS$.
\ede
The aim of this note
is to describe the isomorphism classes of all amalgams that have the
same type as $\cG_{0}$. Note that the classification of amalgams of type $\cG_0$ essentially comes down to classifying all collections of connecting maps $\{\varphi^\sigma_\tau\mid \sigma\sbe \tau\in \sS\}$. Since we classify the amalgams up to isomorphism, some of the $\varphi^\sigma_\tau$ one can specify in advance.

\bde\label{dfn:normalized amalgam} We call an amalgam
$\cG=\{G_\sigma,\varphi^\sigma_\tau\}$ of type $\cG_0$ {\dfn
normalized} if for any simplex $\tau$ we have
$\varphi^{\bar{\tau}}_\tau=\psi^{\bar{\tau}}_\tau$, where
$\bar{\tau}$ is the least maximal face in $\partial \tau$. \ede

\ble\label{lem:many psi's}
Suppose that $\cG$ is normalized as in Definition~\ref{dfn:normalized amalgam}. Suppose that $\sigma\sbe \tau$ have the same largest element $i$.
Then, $\varphi^\sigma_\tau=\psi^\sigma_\tau$.
\ele
\pf
Note that, by induction on $|\sigma|$ and $|\tau|$ we have
 $\varphi^i_\sigma=\psi^i_\sigma$ and $\varphi^i_\tau=\psi^i_\tau$.
Indeed $|\sigma|-1$ applications of the map $\rho\mapsto\bar{\rho}$
 leave the simplex $\{i\}$.
Now since $\psi^i_\tau=\psi^i_\sigma\after\varphi^\sigma_\tau
=\psi^i_\sigma\after\psi^\sigma_\tau$ and all maps are injective, the claim follows.
\qed
\bpr\label{prop:amalgam normalization} Every amalgam of type $\cG_0$
is isomorphic to a normalized amalgam. \epr \pf Let $\cG^{(1)}=\{G_\bullet, \varphi_\bullet^\bullet\}$ be
an arbitrary amalgam of type $\cG_0$. We will construct a normalized amalgam
$\cG^{(2)} =\{G_\bullet, \varsigma_\bullet^\bullet\}$ along with an isomorphism $\phi\colon \cG^{(1)}\to
\cG^{(2)}$. We will define $\varsigma^\sigma_\tau$ and
$\phi_{\tau}$ by induction on the rank of $\tau$. To start the
induction let $\phi_{\tau}=\id$ for all simplices $\tau$ of rank
$0$. Assume that all $\varsigma^\sigma_\tau$ and $\phi_\tau$
have been defined for $\tau$ of rank at most $s\ge 0$. Now let
$\tau$ be a simplex of rank $s+1$. Define
$$\phi_\tau=(\psi^{\bar{\tau}}_\tau)^{-1}\phi_{\bar{\tau}}\ \varphi^{\bar{\tau}}_\tau,$$
 where $\bar{\tau}$ is the least maximal face in $\partial \tau$.

Next, for each $\sigma$ define $\varsigma^\sigma_\tau$ via
Equation~(\ref{eqn:amhom}) to be
$\varsigma^\sigma_\tau=\phi_\sigma\after\varphi^\sigma_\tau\after\phi_\tau^{-1}$.
A direct verification shows that, for any triple $\sigma\sbe\rho\sbe
\tau$, we have
$\varsigma^\sigma_\rho\after\varsigma^\rho_\tau=\varsigma^\sigma_\tau$.

It now follows by definition that $\cG^{(2)}$ is normalized and that
$\phi$ is an isomorphism. \qed

\mn
Proposition~\ref{prop:amalgam normalization} says that we only need to classify normalized amalgams up to
isomorphism.

\bex\label{ex:parabolic amalgam}
Consider a group $G$ acting flag-transitively on a geometry
$\Gamma=(\cO,I,\tau,\ast)$, where $\cO$ denotes the set of objects,
$I$ denotes the set of types, $\typ\colon \cO\to I$ is a type map
and $\ast$ is a symmetric reflexive relation on $\cO$ called the
incidence relation. In the terminology of Buekenhout we shall assume
that $\Gamma$ is connected, transversal, and residually connected.
Let $\xX$ be the simplicial complex in which $V=I$ and
$\sS=\cP(V)-\{\emptyset\}$. Fix a chamber (maximal flag)
$F=(o_i)_{i\in I}$, and for each non-empty subset $\sigma\sbe I$, let
$F_\sigma=(o_j)_{j\in \sigma}$. We now define an amalgam
$\cG_0=\{G_\bullet,\psi_\bullet^\bullet\}$ over $\xX$ setting $G_\sigma =
\Stab_G(F_\sigma)$, for each $\sigma\in \sS$ and letting $\psi^\sigma_\tau\colon
G_\tau\to G_\sigma$ be the inclusion map of subgroups of $G$ whenever $\sigma\sbe \tau$. The group $G_\sigma$ is called the
standard parabolic subgroup of type $\sigma$. A result due to Soul\'e,
Tits, and Pasini now says that $G$ is the universal completion of
$\cG_0$ if and only if the complex whose simplices are the flags of
$\Gamma$ is simply-connected. \eex

\section{Coefficient systems and $1$-cohomology}\label{section:cohomology}
Let $\xX=(V,\sS)$ be a simplicial complex. For any $k\in \NN$, let
$\sS_{k}$ be the set of all simplices of rank $k$ and, for $l\in \NN$,  let $\sS_{\le l}=\bigcup_{0\le k\le l}\sS_k$.

\bde\label{dfn:coefficient system} A {\em coefficient system} on the
simplicial complex $\xX=(V,\sS)$ is a collection
$$\cA=\{A_\sigma, \alpha^\sigma_\tau\mid \sigma\sbe \tau\mbox{ with }\sigma,\tau\in \sS \},$$
where $A_\sigma$ is a group and $\alpha^\sigma_\tau\colon
A_\sigma\to A_\tau$ is a group homomorphism such that whenever
$\sigma\sbe\rho\sbe\tau$, we have
$\alpha^\sigma_\tau=\alpha^\rho_\tau\after\alpha^\sigma_\rho$.
As for amalgams, we shall use the shorthand notation
$\cA=\{A_\bullet,\alpha_\bullet^\bullet\}$. 

\ede
\bde\label{dfn:coeffhom}
We define a {\dfn homomorphism} between coefficient systems
$\cA^{(2)}=\{A_\bullet^{(2)}, {}^{(2)}\alpha_\bullet^\bullet\}$ and
$\cA^{(1)}=\{A_\bullet^{(1)},  {}^{(1)}\alpha_\bullet^\bullet\}$
 to be a map
$\chi=\{ \chi_\sigma\mid \sigma\in\sS\}$ where $\chi_\sigma\colon
A_\sigma^{(2)}\to  A^{(1)}_{\sigma}$ are group homomorphisms  such
that \beq\label{eqn:coeffhom} \chi_{\tau}\after
{}^{(2)}\alpha^\sigma_\tau ={}^{(1)}\alpha^\sigma_\tau\after
\chi_{\sigma}. \eeq We call $\chi$ an {\dfn isomorphism} of coefficient systems if $\chi_\sigma$ is bijective for all $\sigma\in \sS$. \ede

\bde\label{dfn:coeff trivial extension functor}
Adopting the notation from Definition~\ref{dfn:coeffhom}, suppose $\tilde{X}=(\tilde{V},\tilde{\Sigma})$ is a simplicial complex such that $V\sbe \tilde{V}$ and $\Sigma\sbe \tilde{\Sigma}$.
Given a coefficient system $(\cA_\bullet,\alpha_\bullet^\bullet)$ 
 over $X$, we define a coefficient system $(\tilde{\cA},\tilde{\alpha}_\bullet^\bullet)$ over $\tilde{X}$ as follows.
 $$
 \begin{array}{lcl}
 \tilde{A}_\sigma = \left\{\begin{array}{ll} A_\sigma &\mbox{ if }\sigma\in \Sigma\\
 \{1\} & \mbox{ else }\end{array}\right.
 & \mbox{ and }
 \tilde{\alpha}_\tau^\sigma  = \left\{\begin{array}{ll} \alpha_\tau^\sigma & \mbox{ if }\sigma,\tau\in \Sigma\\
 \tilde{A}_\sigma\to \tilde{A}_\tau=\{1\} & \mbox{ else. }
\end{array}\right.
\end{array}
 $$
Now if $\phi=\{\phi_\sigma\mid \sigma\in \Sigma\}\colon \cA^{(2)}\to\cA^{(1)}$ is a homomorphism, then we define $\tilde{\phi}=\{\tilde{\phi}_\sigma\mid \sigma\in \tilde{\Sigma}\}\colon \tilde{\cA}^{(2)}\to\tilde{\cA}^{(1)}$ as follows.
 $$
 \tilde{\phi}_\sigma  = \left\{\begin{array}{ll} 
 \phi_\sigma & \mbox{ if }\sigma \in \Sigma\\
 \id\colon \tilde{A}_\sigma^{(2)}=\{1\}\to \{1\}=\tilde{A}_\sigma^{(1)} & \mbox{ else. }
\end{array}\right.
$$
\ede
\ble\label{lem:coeff trivial extension functor}
With the notation of Definition~\ref{dfn:coeff trivial extension functor} the assignment $\cA\mapsto \tilde{\cA}$ and $\phi\mapsto \tilde{\phi}$ is an embedding of the category of simplicial amalgams over $X$ into the category of coefficient systems over $\tilde{X}$.
In particular, $\phi\colon\cA^{(2)}\stackrel{\cong}{\longrightarrow} \cA^{(1)}$ if and only if $\tilde{\phi}\colon \tilde{\cA}^{(2)}\stackrel{\cong}{\longrightarrow}\tilde{\cA}^{(1)}$.
\ele
\pf
This is a completely straightforward excercise.
\qed

Lemma~\ref{lem:coeff trivial extension functor}
 allows us to replace $X$ by a simplicial complex of rank at least $2$, if necessary. It also allows us to assume that all subsets of $V$ of cardinality $\le k$ are simplices in $\Sigma$.
Given a coefficient system $\cA=\{A_\bullet,\alpha_\bullet^\bullet\}$ on $X$, we
define a cochain complex of pointed sets
$$\cC\colon \cC^0\stackrel{d_0}{\longrightarrow}\cC^1 \stackrel{d_1}{\longrightarrow}\cC^2 $$
where $(\cC^i,\id_i)=\prod_{\sigma\in \sS_i}
(A_\sigma,\id_{A_\sigma})$ as a product of pointed sets,
 $$d_0((a_1,\ldots,a_n)) =(b_{ij}\mid \{i,j\}\in \sS, i<j),\mbox{ where }
b_{ij}=\alpha^{j}_{ij}(a_j^{-1})\alpha^{i}_{ij}(a_i).\\
$$
and
 $$\begin{array}{ll}
   &d_1((a_{ij}\mid \{i,j\}\in \sS, i<j))=(b_{ijk}\mid \{i,j,k\}\in \sS, i<j<k)\\
   \mbox{ where} &\\
 & b_{ijk}  = \alpha^{jk}_{ijk}(a_{jk}^{-1})
  \alpha^{ik}_{ijk}(a_{ik}) \alpha^{ij}_{
  ijk}(a_{ij}^{-1}).
\end{array}$$
Note that the maps $d_i$ are not necessarily group homomorphisms,
although they can be, for instance when the $A_\sigma$'s are
abelian groups. Therefore this is not a chain complex of groups, but
merely a chain complex of pointed sets, where the pointing identifies the identity as a base point in each group $A_\sigma$. It is easy to see that the
maps $d_i$  preserve the base point. 
\ble 
We have $d_1\after d_0 (\cC^0)=(\id_\tau)_{\tau\in \Sigma_2}=\id_{2}$. 
\ele 
\pf It suffices to prove that for any $\sigma$
of rank $0$ and any $\tau$ of rank $2$,  the composition
$$A_\sigma\into \cC^0\stackrel{d_0}{\longrightarrow}\cC^{1} \stackrel{d_{1}}{\longrightarrow}\cC^{2}\to A_\tau$$
sends $A_\sigma$ to $\id_{A_\tau}$. Since the diagram of  $\alpha$'s
is commutative, this is a computation in ordinary cohomology theory
with $\alpha^\sigma_\tau$ applied. \qed

\bde\label{dfn:cocycles} 
For $i=0,1$, the  set of $i$-cocycles of $\cC$ is 
$Z^i(X, \cA)=\{z\in\cC^i\mid d_i(z)=\id_{i+1}\}$. \ede \mn We now define a
right action of $\cC^0$ on $\cC^1$ as follows:
 \beq\label{eqn:coboundary action}
  (b_{ij}\mid 1\le i<j\le n)^{(a_k\mid 1\le k\le n)} =
 (\alpha^{j}_{ij}(a_j^{-1})\cdot b_{ij}\cdot
 \alpha^{i}_{ij}(a_i)\mid 1\le i<j\le n),
  \eeq
   for any
  $(b_{ij}\mid 1\le i<j\le n)\in \cC^1$ and $(a_k\mid 1\le k\le
 n)\in \cC^0$. 
 
\ble\label{lem:d_0}
Let $a,b\in \cC^0$.
Then, we have
\begin{itemize}
\AI{a} $d_0(a)=\id_1^a$, and
\AI{b} $(d_0(a))^b=d_0(ab)$.
\end{itemize}
\ele
\pf
Part $a$ is immediate from the definition, and part (b) follows from part (a) together with the fact that $\cC^0$ acts on $\cC^1$ from the right.
\qed

 \ble\label{lem:C0 preserves Z1} The action of $\cC^0$ on $\cC^1$ preserves
$Z^1(X,\cA)$.
\ele
\pf
Let $a=(a_1,\ldots,a_n)\in \cC_0$ and let
 $z=(z_{ij}\mid 1\le i<j\le n)\in Z^1(X,\cA)$.
 Then
 $z^a = y=(y_{ij}=\alpha^{j}_{ij}(a_j^{-1})\cdot z_{ij}\cdot
 \alpha^{i}_{ij}(a_i)\mid 1\le i<j\le n)$.
The projection of the co-boundary $d_1(z^a)$ on $A_\tau$ for some  simplex $\tau=\{i,j,k\}$ with $i<j<k$ equals:
$$
 \alpha^{jk}_{ijk}(y_{jk}^{-1})
  \alpha^{ik}_{ijk}(y_{ik}) \alpha^{ij}_{
  ijk}(y_{ij}^{-1})\\
 $$
 Using the definition of $y$ we get
 $$\begin{array}{rcl}
 \alpha^{jk}_{ijk}\left(
 \alpha^{j}_{jk}(a_j^{-1})
\cdot z_{jk}^{-1}
\cdot
\alpha^{k}_{jk}(a_k)
 \right)
 &\cdot &
  \alpha^{ik}_{ijk}\left(
  \alpha^{k}_{ik}(a_k^{-1})\cdot z_{ik}\cdot
 \alpha^{i}_{ik}(a_i)
  \right) \\
  & \cdot &
\alpha^{ij}_{ijk}\left(
\alpha^{i}_{ij}(a_i^{-1})
\cdot z_{ij}^{-1}
\cdot
\alpha^{j}_{ij}(a_j)
  \right)
    \end{array}\\
$$
We use the composition property of the $\alpha$-homomorphisms and some cancellations to get:
$$
\alpha^{j}_{ijk}(a_j^{-1})
\cdot 
  \alpha^{jk}_{ ijk}(
z_{jk}^{-1})
\cdot
    \alpha^{ik}_{ijk}(
z_{ik})
\cdot  \alpha^{ij}_{ijk}(
z_{ij}^{-1})
\cdot
 \alpha^{j}_{ijk}(a_j),  \\
$$
and this equals $\id_{A_{ijk}}$ since $z$ is a cocycle.
\qed

\bde\label{dfn:cohomology} 
The zero cohomology set is $H^0(X,\cA)=Z^0(X,\cA)$.
The orbits in
  $Z^1(X,\cA)$ under the action defined above are called
  $1$-cohomology classes. The first cohomology set is the collection of $\cC_0$ orbits of $1$-cocycles. We write
 $$H^1(X,\cA)=Z^1(X,\cA)^{\cC^0}.$$
Moreover, for each $z\in Z^i(X,\cA)$ we denote its cohomology class as $[z]$. Note that for $z\in Z^0(X,\cA)$ we have $[z]=\{z\}$.
\ede

\bco
The zero cohomology set $H^0(X,\cA)$ is a group.
\eco
\pf
This follows from Lemma~\ref{lem:d_0} since $d_0(ab)=((\id_1)^a)^b=\id_1^b=\id_1$.
\qed
\subsection{A cohomology exact sequence}
\bde\label{dfn:normal subsystem} A normal subsystem of $\cA$ is a collection $\{N_\sigma\mid \sigma\in \sS\}$ of normal subgroups
 $N_\sigma\normal A_\sigma$ such that $\cN=\{N_\bullet,\alpha_\bullet^\bullet\}$ is a coefficient system on $X$. We shall denote this as
  $\cN\normal \cA$.
The normal sub coefficient system  naturally gives rise to a
quotient coefficient system by setting
$\cA/\cN=\{A_\bullet/N_\bullet, \alpha_\bullet^\bullet\}$.
Note that for $i=0,1$, $\cC^i(\cN)\normal \cC^i(\cA)$ and therefore,  $\cC^i(\cA/\cN)\cong\cC^i(\cA)/\cC^i(\cN)$.
\ede

\bth\label{thm:H0 on H1} 
The group $H^0(X,\cA/\cN)$ acts on $H^1(X,\cN)$ in a natural way. The orbits of this action are precisely the fibers of the map $H^1(X,\cN)\stackrel{i_1}\to H^1(X,\cA)$, which takes each $\cC^0(\cN)$ orbit on $Z^1(X,\cN)$ to the unique $\cC^0$ orbit on $Z^1(X,\cA)$ that contains it.
\eth
\pf
Given $\bar{a}=(a_\sigma N_\sigma)_{\sigma\in\sS_0}\in H^0(X,\cA/\cN)$ and 
 $n=(n_\tau)_{\tau\in \sS_1}\in Z^1(X,\cN)$, let $a=(a_\sigma)_{\sigma\in\sS_0}\in \cC^0(\cA)$. We define
  $$[n]^{\bar{a}}=[n^a].$$
We claim that this is well-defined. 
Indeed, let $m\in \cC^0(\cN)$ and $a\in \cC^0(\cA)$ and let $m'\in \cC^0(\cN)$ be such that $am=m'a$.
Then, $n^{am}=n^{m'a}=(n^{m'})^a$ and so
 $[n^{am}]=[(n^{m'})^a]=[n^a]$.
Suppose now that moreover, $\bar{a}\in Z^0(X,\cA/\cN)$. Then, 
 $d_0(a)\in \cC^1(\cN)$.
This means that for each $\{i,j\}\in \sS_1$ with $i<j$, we have
 $m_{ij}=(d_0(a))_{ij}=\alpha^j_{ij}(a_j^{-1})\alpha^{i}_{ij}(a_i)\in N_{ij}$. 
Hence, if $n=(n_{ij})_{\{i,j\}\in\sS_1}$, then, 
\beq\label{eqn:mij}
 n^a=(\alpha^j_{ij}(a_j^{-1})\ n_{ij}\ \alpha^i_{ij}(a_i))_{\{i,j\}\in \sS_1}
 =(\alpha^j_{ij}(a_j^{-1})\ n_{ij}\ \alpha^j_{ij}(a_j) m_{ij})_{\{i,j\}\in \sS_1}\eeq
 and this belongs to $\cC^1(\cN)$ since $N_{ij}\normal A_{ij}$.
This concludes the proof of our claim.

By definition $i_1([n^a])$ and $i_1([n])$ are in the same cohomology class of $H^1(X,\cA)$. Conversely, suppose that
 $[n]$ and $[n']$ are in $H^1(X,\cN)$ such that $i_1([n])=i_1([n'])$.
 Then there is some $a\in \cC^0(X,\cA)$ with 
  $n'=n^a$. We claim that $d_0(a)\in \cC^1(\cN)$ so that $[n]$ and $[n']$ are in the same $H^0(X,\cA/\cN)$-orbit. 
Indeed define $d_0(a)=(m_{ij})_{\{i,j\}\in \sS_1}$.
Then Equation~(\ref{eqn:mij}) still holds for $[n']=[n^a]$.
Since  $n'$ and $n$ belong to $\cC^1(\cN)$ and 
$\cC^1(\cN)\normal \cC^1(\cA)$, also $m\in \cC^1(\cN)$.
\qed

\bth\label{thm:exact sequence}
For any $\cN\normal \cA$ there is a natural exact sequence of pointed sets
$$
0\to 
H^0(X,\cN)\stackrel{i_0}{\to} H^0(X,\cA)\stackrel{\kappa_0}{\to} H^0(X,\cA/\cN)\stackrel{\delta^*}{\to} 
H^1(X,\cN)\stackrel{i_1}{\to} H^1(X,\cA)\stackrel{\kappa_1}{\to} H^1(X,\cA/\cN) 
$$
\eth
\pf
For $j=0,1$, the map $i_j$ is given by the inclusion maps
 $N_\sigma\into A_\sigma$ and the map $\kappa_j$
  is given by the canonical homomorphism $A_\sigma\to A_\sigma/N_\sigma$ for any $\sigma\in \sS$.
The map $\delta^*$ is defined as 
 $$\delta^*((a_\sigma N_\sigma)_{\sigma\in \sS_0})=d_0((a_\sigma)_{\sigma\in\sS_0}).$$
That this is well-defined can be seen as follows.
Since the $\alpha^\bullet_\bullet$ are group homomorphisms that preserve $\cN$, it  follows that the following diagram commutes:
$$\xymatrix{
0 \ar[r] & \cC^0(\cN) \ar[r]\ar[d]^{d_0} & \cC^0(\cA)\ar[r]\ar[d]^{d_0} & \cC^0(\cA/\cN)\ar[r]\ar[d]^{d_0} & 0\\
0 \ar[r] & \cC^1(\cN) \ar[r] \ar[d]^{d_1} & \cC^1(\cA)\ar[r]\ar[d]^{d_1} & \cC^1(\cA/\cN)\ar[r]\ar[d]^{d_1} & 0\\
0 \ar[r] & \cC^2(\cN) \ar[r] & \cC^2(\cA)\ar[r] & \cC^2(\cA/\cN)\ar[r] & 0\\
}$$
Note that the $d_0$'s and $d_1$'s are merely maps of pointed sets. Note that $d_1\after d_0=0$ and that the rows in the diagram are exact sequences of group homomorphisms.
A pointed-set version of the Snake Lemma shows that $\delta^*$ is well-defined and that the sequence is exact up to $H^1(X,\cN)$.
Exactness at $H^1(X,\cN)$ follows from Theorem~\ref{thm:H0 on H1}.
It is easy to see that $\im(i_1)\le \ker(\kappa_1)$.
Let $[a]\in\ker(\kappa_1)$.
Since $a\in Z^1(X,\cA)$, we have $d_1(a)=\id_2$, and since $[a]\in \ker(\kappa_1)$ there exists some $b\in \cC^0(\cA)$ such that
 $a^b=n\in \cC^1(\cN)$.
By Lemma~\ref{lem:C0 preserves Z1} we know that
 $d_1(n)=\id_2$ and so $[a]=[n]\in H^1(X,\cN)$.
\qed

\ble\label{lem:quotient
cochain complex} 
Assume $X$ is a $2$-simplex and $\cN\normal \cA$ is a normal subsystem such that, for each $\sigma\sbe\tau\in\sS$ $\alpha^\sigma_\tau\colon N_\sigma\to N_\tau$ is an isomorphism. Then
$H^1(X, \cA)=H^1(X, \cA/\cN)$. 
\ele
\pf
By Theorem~\ref{thm:exact sequence} it suffices to show that $H^1(X,\cN)=0$ and that $\kappa_1$ is onto.
Without loss of generality assume that $X$ is the set of non-empty subsets of $\{1,2,3\}$.
Let $n=(n_{ij})\in Z^1(X,\cN)$. That means that $\alpha^{23}_{123}(n_{23}^{-1})\alpha^{13}_{123}(n_{13})\alpha^{12}_{123}(n_{12}^{-1})=0$.
Because of this, and since all $\alpha$ maps are isomorphisms we can find $m=(m_1, m_2,m_3)$ such that  
$$\begin{array}{rl}
\alpha^3_{23}(m_3^{-1})\alpha^2_{23}(m_2)&=n_{23},\\
\alpha^3_{13}(m_3^{-1})\alpha^1_{13}(m_1)&=n_{13},\\\alpha^2_{12}(m_2^{-1})\alpha^1_{12}(m_1)&=n_{12}.\\
\end{array}$$
Then, $d_0(m)=\id_1^{m}=n$. Thus $H^1(X,\cN)=0$.
Now take $\bar{a}=(\bar{a_{ij}})\in Z^1(X,\cA/\cN)$. This means that 
 for a representative $a=(a_{ij})$ we have  $d_1(a)=n\in N_{123}$.
Now let  $a'=(a'_{ij})$ be given by 
 $a'_{12}=a_{12} m^{a_{12}}$, where $n=\alpha^{12}_{123}(m)$ and $a'_{ij}=a_{ij}$ otherwise.
Now 
$$d_1(a')
= \alpha^{23}_{123}(a_{23}^{-1})\alpha^{13}_{123}(a_{13})\alpha^{12}_{123}(a_{12}^{-1})n^{-1}=\id_{123}.$$
Clearly $(\bar{a_{ij}'})=(\bar{a_{ij}})$ so we are done.
\qed

\subsection{The coefficient system of an amalgam}\label{sec:the coefficient system of an amalgam}
We now construct a coefficient system
$\cA_0=\{A_\bullet,\alpha_\bullet^\bullet\}$, for the amalgam $\cG_0$ by
setting for each $\sigma\in \sS$
 $$A_\sigma=\{g\in \Aut(G_\sigma)\mid g(\bar{G}_{\sigma,\tau})=\bar{G}_{\sigma,\tau}\mbox{ for all }\tau\in\sS\mbox{ with }\sigma\sbe\tau\}$$
(in particular, if $\sigma$ is maximal, this implies $A_\sigma=\Aut(G_\sigma)$).
Moreover, for each pair $(\sigma,\tau)$ with $\sigma\sbe\tau$ we
define $\alpha^\sigma_\tau \colon A_\sigma\to A_\tau$ given by
$\ad(\psi^\sigma_\tau)$, where $\ad(x)(y)=x^{-1}yx$. 
If $\phi=\{\phi_\sigma\mid \sigma\in\Sigma\}\colon \cG^{(1)}\to \cG^{(2)}$ is a homomorphism of simplicial amalgams, then 
 $\chi=\{\chi_\sigma=\ad(\phi_\sigma)\mid \sigma\in \Sigma\}\colon \cA_0^{(2)}\to\cA_0^{(1)}$ is a homomorphism of coefficient systems. One verifies easily that this assignment defines a contravariant functor from the category of simplicial amalgams over $X$ to the category of coefficient systems over $X$.
 
We now let $\cC_0^\bullet$ be the cochain
complex associated to $\cA_0$.

\bex\label{ex:GL42 1-cohomology} In Example~\ref{ex:parabolic
amalgam} take $\Gamma$ to be the projective $3$-space $\PG(V)$,
where $V$ has dimension $4$ over $\FF_2$ and let $G=\SL_4(2)$. Thus
$\cG_0$ is an amalgam over a complex $\xX=(V,\sS)$, where $\sS$
consists of all non-empty subsets of $V=\{1,2,3\}$. A computation
with GAP now reveals the following: $\alpha^{1}_{12}\colon
A_1 \to A_{12}$  and $\alpha^{3}_{23}\colon A_3 \to
A_{23}$ are surjective. This implies that every element in
$\cC_0^1$ is in a cohomology class with an element
$(\id_{12},\id_{23},a_{13})$, for some
$a_{13}\in A_{13}$. However,  another calculation shows
that $\alpha^{13}_{123}\colon A_{13} \to
A_{123}$ is an isomorphism. Now for $(\id_{12},\id_{23},a_{13})\in Z^1(X,\cA_0)$ the $1$-cocycle condition
forces $a_{13}=\id_{13}$ so that $H^1(X,\cA_0)=0$.
\eex
\section{The correspondence between cohomology classes and amalgams}\label{section:cohomology and amalgams}
Consider the reference amalgam
$\cG_0=\{G_\bullet, \psi_\bullet^\bullet\}$ over the connected simplicial complex $X=(V,\Sigma)$.
Using Lemma~\ref{lem:trivial extension functor} we can assume that $X$ has rank at least $2$ and contains all $3$-subsets of $V$ as a simplex.
Recall that every amalgam of
type $\cG_0$ has the same target subgroups $\bar{G}_{\sigma,\tau}$.
Consider now $\cG=\{G_\bullet, \varphi_\bullet^\bullet\}$ a normalized
amalgam of type $\cG_0$. 
For each $\sigma =\{i,j\}
\in \Sigma_1$ with $i<j$ define $a^{\cG}_{ij}=
(\varphi^{i}_{ij})^{-1}\after \psi^{i}_{ij}$. Note that
it follows that $a_\sigma \in A_\sigma$ and so the collection
$a=\{a^{\cG}_\sigma\mid \sigma \in \Sigma_1\}$ is an element of
$\cC^1$. \bpr\label{prop:amalgams to cycles} The collection
$a=\{a^{\cG}_\sigma\mid \sigma \in \Sigma_2\}$ is an element of
$Z^1(X,\cA_0)$. Moreover the correspondence $\cG \to
\{a^\cG_\bullet\}$ is a bijection between the set of  normalized
amalgams of type $\cG_0$ and the set $Z^1(X,\cA_0)$ \epr

\pf
We first need to prove that $a\in Z^1(X,\cA_0)$. Let us consider
$\tau=\{i_1,i_2, i_3\}\in \Sigma_2$,  such that $i_1<i_2< i_3$. In
order to have a cleaner set of notations (and no double subscripts)
we will assume that $i_1,i_2,i_3= 1,2,3$. The general case is
similar. The amalgam $\cG$ is normalized so
$\varphi^{23}_{123}=\psi^{23}_{123},
\varphi^{2}_{12}=\psi^{2}_{12},
\varphi^{3}_{23}=\psi^{3}_{23}$ and
$\varphi^{3}_{13}=\psi^{3}_{13}$. We also have that
$\varphi^{1}_{12}=\psi^{1}_{12}a_{12}^{-1},
\varphi^{1}_{13}=\psi^{1}_{13}a_{13}^{-1}$ and
$\varphi^{2}_{23}=\psi^{2}_{23}a_{23}^{-1}$.

The information above is summarised in the following diagram
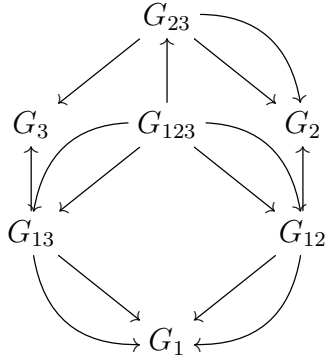
\[\xymatrix{
& G_{23} \ar[dl] \ar[dr]\ar@/^1.5pc/[dr]&\\
 G_{3} & G_{123}\ar[dl]\ar@/_1.5pc/[dl]\ar[u] \ar@/^1.5pc/[dr] \ar[dr]& G_{2}\\
 G_{13}\ar[u]  \ar[dr] \ar@/_1.5pc/[dr]&  & G_{12}\ar@/^1.5pc/[dl]\ar[dl]\ar[u]\\
&G_{1}& 
}\]
 
  Here the straight arrows signify the maps $\psi^\tau_\sigma$ and the curved maps signify the maps $\varphi^\tau_\sigma$, wherever they differ. 
Thus, the diagram obtained by considering only straight arrows  corresponds to the amalgam $\cG_0$ and the diagram obtained by taking the curved arrows wherever possible corresponds to the amalgam $\cG$. Both these diagrams commute. 
Note that by the
 commutativity  of the squares, a walk in the diagram for $\cG$ from
 $G_{123}$ to $G_{123}$ along the path $$\{1,2,3\}\to
 \{1,2\} \to \{1\} \to \{1,3\} \to \{3\} \to \{2,3\} \to \{2\} \to
 \{1,2\} \to \{1,2,3\}$$ gives
$$\id_{G_{123}}=(\psi^{12}_{123})^{-1}(\vph^{2}_{12})^{-1}\vph^{2}_{23} (\vph^{3}_{23})^{-1} \vph^{3}_{13} (\vph^{1}_{13})^{-1} \vph^{1}_{12}\psi^{12}_{123}.$$
Using, that, for $i<j$, we have $\varphi^i_{ij}=\psi^i_{ij}a_{ij}^{-1}$
 and $\varphi^j_{ij}=\psi^j_{ij}$,  we have
$$\begin{array}{ll}
\id_{G_{123}}& =(\psi^{12}_{123})^{-1}(\psi^{2}_{12})^{-1}(\psi^{2}_{23}a_{23}^{-1})(\psi^{3}_{23})^{-1}(\psi^{3}_{13})(a_{13}(\psi^{1}_{13})^{-1})(\psi^{1}_{12}a_{12}^{-1})\psi^{12}_{123}.\\
\end{array}$$
Then, using that the diagram commutes we find
$$\begin{array}{ll}
   \id_{G_{123}} & =(\psi^{12}_{123})^{-1}\psi^{12}_{123}(\psi^{23}_{123})^{-1}a_{23}^{-1}\psi^{23}_{123}(\psi^{13}_{123})^{-1}a_{13}\psi^{13}_{123}(\psi^{12}_{123})^{-1}a_{12}^{-1}\psi^{12}_{123}\\ 
    & \\
    & =(\psi^{23}_{123})^{-1}a_{23}^{-1}\psi^{23}_{123}(\psi^{13}_{123})^{-1}a_{13}\psi^{13}_{123}(\psi^{12}_{123})^{-1}a_{12}^{-1}\psi^{12}_{123}\\
    
    & \\
    & =\alpha^{23}_{123}(a_{23}^{-1})\alpha^{13}_{123}(a_{13})\alpha^{12}_{123}(a_{12}^{-1}).\\
\end{array}
$$
This last equation means that for any amalgam $\bar{G}_\bullet$, the collection $a^\cG_\bullet$ is a 1-cocycle.

Conversely suppose $a=\{a_\sigma \mid \sigma \in \Sigma_1\}$ is a
1-cocycle. We need to construct an amalgam $\cG$ so that
$a=\{a^\cG_\bullet\}$.
We first define $\varphi_\tau^\sigma$ for any $\sigma\sbe \tau\in\sS$. Let $\max\sigma=j$ and $k=\max \tau$.
It follows from Lemma~\ref{lem:many psi's} that we must define
  $$\varphi_\tau^\sigma=\psi_\tau^\sigma\mbox{ if }j=k.$$
We also define
  $$\varphi_\tau^\sigma=(\psi_\sigma^j)^{-1}\after\psi_{jk}^j a_{jk}^{-1}\after \psi_\tau^{jk}\mbox{ if }j<k.$$
Note that this definition is also forced upon us. 
Namely, for $\sigma=\{j\}$ and $\tau=\{j,k\}$ this is the only possibility since we insist that $a=\{a^\cG_\bullet\}$.
Moreover, normality already forced us to set $\varphi_\sigma^j=\psi_\sigma^j$
and $\varphi_\tau^{\{j,k\}}=\psi_\tau^{\{j,k\}}$. Hence our definitions of $\varphi_\tau^j$ and $\varphi_\tau^\sigma$ are forced upon us by the requirements
 $$\begin{array}{ccc}
 \varphi_\tau^j=\varphi^j_{\{j,k\}}\after\varphi_\tau^{\{j,k\}}=\varphi_\sigma^j\after\varphi_\tau^\sigma.
 \end{array}
 $$
It now suffices to show that for any $\rho$ with $\rho\sbe\sigma\sbe \tau$ we have
\beq\label{eqn:consistent} \varphi_\tau^\rho=\varphi^\rho_\sigma\after \varphi^\sigma_\tau.\eeq
Let $i=\max \rho$ and let $j$ and $k$ be as above.
If $i=j=k$ then (\ref{eqn:consistent}) follows since all $\psi$ commute.
If $\{i,j,k\}=\{i,k\}$, then either $\varphi_\tau^\sigma=\psi_\tau^\sigma$ or $\varphi_\sigma^\rho=\psi_\sigma^\rho$.
Suppose the latter holds.
Then, 
$$\varphi_\tau^\rho=
(\psi^i_\rho)^{-1}\after\psi^i_{ik}a_{ik}^{-1}\after \psi_\tau^{ik}=
\psi^\rho_\sigma (\psi^i_\sigma)^{-1}\after\psi^i_{ik}a_{ik}^{-1}\after \psi_\tau^{ik}=\varphi_\sigma^\rho\after\varphi^\sigma_\tau$$
and a similar argument holds in the former case.
Finally let $i<j<k$.
Then, (\ref{eqn:consistent}) amounts to 
$$
(\psi^i_\rho)^{-1}\after\psi^i_{ik}a_{ik}^{-1}\after \psi_\tau^{ik}=
(\psi^i_\rho)^{-1}\after\psi^i_{ij}a_{ij}^{-1}\after \psi_\sigma^{ij}
(\psi^j_\sigma)^{-1}\after\psi^j_{jk}a_{jk}^{-1}\after \psi_\tau^{jk}$$
Multiplying by $(\psi^i_{ijk})^{-1}\psi_\rho^i$ on the left
 and $(\psi_\tau^{ijk})^{-1}$ on the right, replacing 
 $\psi_\sigma^{ij}(\psi_\sigma^j)^{-1}=(\psi^j_{ij})^{-1}$, and using that $(\psi^j_{ij})^{-1}\psi^j_{jk}=\psi^{ij}_{ijk}(\psi^{jk}_{ijk})^{-1}$ and the $\psi$ all commute this reduces to 
  $$
  \alpha^{ik}_{ijk}(a_{ik}^{-1})
  =
\alpha^{ij}_{ijk}(a_{ij}^{-1})
  \alpha^{jk}_{ijk}(a_{jk}^{-1})$$
and this is equivalent to the fact that 
 $a$ is a cocycle.
We already demonstrated that, for any $1$-cocycle $a$ there is (at most) a unique normalized amalgam $\cG$ with $a=\{a^\cG_\bullet\}$, so we are done.
\qed  
\subsection{Isomorphisms and co-boundaries}\label{subsec:coboundaries}
\bpr\label{prop:isomorphism vs cohomologous}
Two normalized amalgams of type $\cG_0$ are isomorphic if and only if the corresponding $1$-cocycles are cohomologous.
\epr
\pf
For $l=1,2$, let ${}^{(l)}\cG=\{G_\bullet,{}^{(l)}\varphi^\bullet_\bullet\}$ be a normalized amalgam of type $\cG_0$ corresponding to
 a cocycle $z^{(l)}=\{z_{\sigma}^{(l)}\mid\sigma\in \sS_1\}$.
Recall that this means that, for $\{i,j\}\in \sS_1$ with $i<j$ we have
 $${}^{(l)}\varphi^{i}_{ij}=\psi^{i}_{ij}(z_{ij}^{(l)})^{-1}.$$

Suppose $\phi:\cG^{(1)}_\bullet \to
\cG^{(2)}_\bullet$ is an isomorphism. It then follows that
${}^{(l)}\varphi_\tau^{\bar \tau}= \psi_\tau^{\bar \tau}$ for all
$\tau$ with $|\tau|>1$. In particular Equation~\ref{eqn:amhom}
becomes
$\phi_{\bar \tau}\after \psi^{\bar \tau}_\tau
=\psi^{\bar \tau}_\tau\after \phi_{\tau}$
and so we have 
\beq\label{eqn:vertex maps determine everything}
\phi_\tau =\alpha^{\bar\tau}_{\tau}(\phi_{\bar \tau})
\eeq
Moreover if $i<j$, we get the following commutative diagram:
\beq\label{eqn:rank 2}
\xymatrix{
 & G_{i}  \ar[r]^{\phi_{i}} & G_{i} & \\
  & G_{ij}\ar[u]^-{\psi^{i}_{ij} \after (z^{(1)}_{ij})^{-1}= ^{(1)}\varphi^i_{ij}}  \ar[r]_{\phi_{ij}\\ =\alpha^j_{ij}(\phi_j)} &   G_{ij} \ar[u] _-{^{(2)}\varphi^i_{ij}=\psi^{i}_{ij}\after (z^{(2)}_{ij})^{-1}} &\\ 
}
\eeq
If we compare the two cocycles $z^{(l)}_{ij}=({}^{(l)}\varphi^{i}_{ij})^{-1}\after\psi^{i}_{ij}$ we see that
$$z^{(1)}_{ij}= ({}^{(1)}\varphi^{i}_{ij})^{-1}\after\psi^{i}_{ij} =\alpha^{j}_{ij}(\phi_j^{-1})\after ({}^{(2)}\varphi^{i}_{ij})^{-1}\after \phi_i \after \psi^{i}_{ij}=\alpha^{j}_{ij}(\phi_j^{-1})\after z^{(2)}_{ij}\after \alpha^{i}_{ij}(\phi_i) .$$
This shows that in fact the diagram~(\ref{eqn:rank 2}) is commutative if and only if $z^{(1)}=(z^{(2)})^{\{\phi_k\mid 1\le k\le n\}}$.
In particular, $z^{(1)}$ and $z^{(2)}$ are cohomologous.
(Note here that since $\phi$ preserves all $G_\sigma$
  and $\bar{G}_\sigma^\tau$, we have $\phi_k\in A_k$, for all $k$.)

\mn
Conversely, suppose that $z^{(1)}$ and $z^{(2)}$ are cohomologous, that is, they belong to the same $\cC^0$-orbit.
Let $f=\{f_k\mid k\in V\}\in \cC^0$ such that $z^{(1)}= (z^{(2)})^f$.
We now define an isomorphism $\phi\colon {}^{(1)}\cG\to {}^{(2)}\cG$.
First set
 $$\phi_{k} = f_k\mbox{ for all }1\le k\le n.$$ 
We now note that since $^{(1)}\cG$ and $^{(2)}\cG$ are normalized of the same type, 
Equation~(\ref{eqn:vertex maps determine everything}) must be satisfied and inductive use together with the composition properties of the $\alpha$ maps, shows that, for each simplex $\tau=\{i_1,\ldots,i_m\}$, with $i_1<\cdots<i_m$, we must have 
$$\phi_\tau=\alpha^{i_m}_\tau(f_{i_m}).$$
It now suffices to check that for all simplices $\sigma\sbs \tau$ we have
 $$\phi_\sigma\after {}^{(1)}\varphi^\sigma_\tau= {}^{(2)}\varphi^\sigma_\tau\after\phi_\tau.$$
For $\sigma=\{i\}$ and $\tau=\{i,j\}$ with $i<j$, this requires the diagram~(\ref{eqn:rank 2}) to be commutative and we already saw that this is equivalent to 
 $z^{(1)}=(z^{(2)})^{\{\varphi_{k}\mid 1\le k\le n\}}$.
Thus for all $\tau$ with $|\tau|=2$, we're done.

Next, consider $\sigma\sbs\tau$ where $|\tau|>2$.
Suppose that $i$ and $j$ are the largest vertices of $\sigma$ and $\tau$ respectively.
If $i=j$, then 
 $$\phi_\tau=\alpha^{i}_{\tau}(f_i)=\alpha^\sigma_\tau \after\alpha^{i}_\sigma(f_i)=\alpha^\sigma_\tau(\phi_\sigma)=
 (\psi^\sigma_\tau)^{-1}\after\phi_\sigma\after \psi^\sigma_\tau=
 (^{(2)}\varphi^\sigma_\tau)^{-1}\after\phi_\sigma\after {}^{(1)}\varphi^\sigma_\tau$$
and we are done. Here the last equality follows from Lemma~\ref{lem:many psi's}.

If $i\ne j$, then $i<j$.
Since, for $l=1,2$,  ${}^{(l)}\cG$ is normalized, we have
 $^{(l)}\varphi^{i}_{\sigma}=\psi^{i}_{\sigma}$ and 
 $^{(l)}\varphi^{ij}_{\tau}=\psi^{ij}_{\tau}$, again by Lemma~\ref{lem:many psi's}.
It follows that the left and right square in the diagram below are commutative.
$$
\xymatrix{
G_\sigma \ar[drr]^{^{(2)}\varphi^{i}_\sigma=\psi^{i}_\sigma} 
&&&&
&& \ar[dll]_-{^{(2)}\varphi^{ij}_\tau=\psi^{ij}_{\tau}}   G_\tau \ar[llllll]^{^{(2)}\varphi^\sigma_\tau}\\
&& G_{i} 
&&\ar[ll]_{^{(2)}\varphi^{i}_{ij}}   G_{ij} 
&&\\
&&G_{i}   \ar[u]^{f_i=\phi_{i}}  
&& \ar[ll]^{^{(1)}\varphi^{i}_{ij}}G_{ij}\ar[u]_{\phi_{ij}=\alpha^{j}_{ij}(f_j)}
&& \\
G_\sigma \ar[urr]_{\hspace{2pt} ^{(1)}\varphi^{i}_\sigma=\psi^{i}_{\sigma}} \ar[uuu]^{\alpha^{i}_\sigma(f_i)=\phi_\sigma} 
&&&&
&& G_\tau  \ar[ull]^-{^{(1)}\varphi^{ij}_\tau=\psi^{ij}_{\tau}\hspace{9pt}}  \ar[uuu]_{\phi_\tau=\alpha^{j}_\tau(f_j)} \ar[llllll]_{^{(1)}\varphi^\sigma_\tau}\\
}
$$
Moreover, the middle square is commutative because of the rank-$2$ case.
The top and bottom square are also commutative and so the result follows.
\qed

\mn
The proof of Theorem~\ref{mainthm1} is now complete. 

\bex\label{ex:GL42 amalgams}
In Example~\ref{ex:GL42 1-cohomology} we considered the special case of Example~\ref{ex:parabolic amalgam}, where the amalgam $\cG_0$ consists of standard parabolic subgroups of $\SL_4(2)$ as it acts on the projective $3$-space $\PG(V)$, where $V=\FF_2^4$.
We computed that $H^1(X,\cA)=0$, so that by Theorem~\ref{mainthm1} there exists a unique amalgam of this type.
\eex
\section{Amalgams over small complexes}\label{applications}
\subsection{Goldschmidt's Lemma}

The simplest case of the theory is the celebrated Goldschmidt's Lemma. It arises as follows from our setup. Let $X$ be the $1$-simplex $\{\{1\},\{2\}, \{1,2\}\}$ and let $\cG_0$ be a generic amalgam over $X$. We denote its three groups by $G_1, G_2, G_{12}$ and let $\psi^{i}\colon G_{12}\to G_{i}$, for $i=1,2$.
We also define the three automorphisms groups as in Subsection~\ref{sec:the coefficient system of an amalgam} as $A_{12}
=\Aut(G_{12})$ and $$A_i=\{g \in \Aut(G_i) |  g(\bar G_{i,\{1,2\}})=\bar G_{i,\{1,2\}}\}$$
 respectively $\bar A_i = \ad(\psi^{i})(A_i)\le A_{12}$.
Theorem~\ref{mainthm1}  now reads as
\bco{\rm (Goldschmidt's Lemma, see~\cite[\S 2.7]{Gol80} 
and~\cite[Ch.~16]{DelGolSte1985})} \label{lem:Goldschmidt} There is a 1-1 correspondence between isomorphism classes of amalgams of type $\cA_0=\{G_1, G_2, G_{12}\}$ and double cosets of $\bar A_{1}$, and $\bar A_{2}$ in $A_{12}$.
\eco

\subsection{Triangular complexes}
Let $X$ be the $2$-dimensional simplex consisting of all non-empty subsets of $V=\{1,2,3\}$ and let $\cG_0$ be a generic amalgam over $X$. Such an amalgam  arises naturally from a group acting flag-transitively on a rank $3$ geometry as in Example~\ref{ex:parabolic amalgam}.

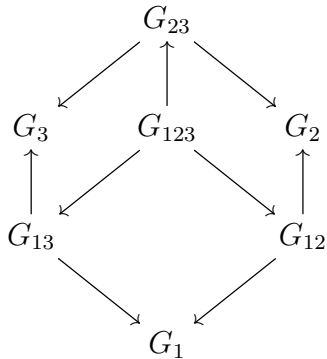
\[\xymatrix{
& G_{23} \ar[dl] \ar[dr]&\\
 G_{3} & G_{123}\ar[dl]\ar[u] \ar[dr]& G_{2}\\
 G_{13}\ar[u]  \ar[dr]&  & G_{12}\ar[dl]\ar[u]\\
&G_{1}&
}\]
In order to apply Theorem~\ref{mainthm1}, we need to consider the corresponding  coefficient system $\cA$.

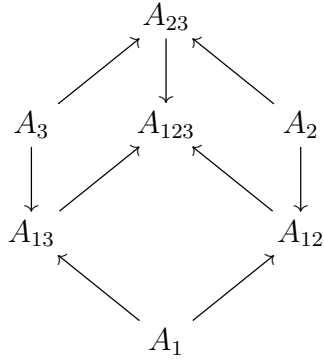
\[\xymatrix{
& A_{23}\ar[d]&\\
 A_{3}\ar[ur]\ar[d] & A_{123}& A_{2}\ar[ul]\ar[d]\\
 A_{13}\ar[ur]&  & \ar[ul]A_{12}\\
&\ar[ul]A_{1}\ar[ur]& }\]

Here $A_{i}=\{g \in \Aut(G_{i}) \mid 
g(\bar{G}_{ij}) = \bar{G}_{ij}\mbox{ for all }j\ne i \}$, $A_{ij}=\{g \in
\Aut(G_{ij}) \mid g(\bar{G}_{123})= \bar{G}_{123}$,
$A_{123}=\Aut(G_{123})$ and the arrows represent the
maps $\alpha^\sigma_\tau\colon A_\sigma\to A_\tau$, which can be
viewed as restriction maps.

We shall consider $N=\Inn(G_{123})$ and note that for any
$\sigma$ there is a group $N_\sigma \lhd A_\sigma$ so that $N\cong
N_\sigma$ and that $\alpha^\sigma_\tau$ restricted to $N_\sigma$
gives an isomorphism between $N_\sigma$ and $N_\tau$.  We define
then $\bar A_\sigma =A_\sigma/N_\sigma$ and
$\bar\alpha^\sigma_\tau\colon \bar A_\sigma\to \bar A_\tau$. Define
then the set

$$\bar{\cC}^1 =\{(a_{23},a_{13},a_{12})\mid a_{ij} \in
\bar A_{ij} \mbox{ and } \bar\alpha^{23}_{
123}(a_{23}^{-1})
  \bar \alpha^{13}_{123}(a_{13}) \bar \alpha^{12}_{
  123}(a_{12}^{-1})=id_{\bar A_{123}}\}$$

  The group $\bar{\cC}^0=\bar A_{1}\times \bar A_{2}\times \bar
A_{3}$ acts on $\bar{\cC}^1$ via

$$\begin{array}{l}(a_{23},a_{13},a_{12})^{(d_1,d_2,d_3)}=\\(\bar \alpha^{3}_{23}(d_3^{-1})a_{23}\bar \alpha^{2}_{23}(d_2),
\bar \alpha^{3}_{13}(d_3^{-1})a_{13}\bar
\alpha^{1}_{23}(d_1), \bar
\alpha^{2}_{12}(d_2^{-1})a_{12}\bar
\alpha^{1}_{12}(d_1)) \end{array}$$

\bpr The isomorphism classes of amalgams of type $\cG_0$ are in
bijection to the orbits of $\bar{\cC}^0$ on $\bar{\cC}^1$ \epr

\pf Immediate from Theorem~\ref{mainthm1} and
Lemma~\ref{lem:quotient cochain complex}\qed

\bre
Note that the group and the action is exactly as in the example of
triangular rank two amalgams from \cite{BloHof2011}, however the
existence of the group $G_{123}$ gives smaller $A's$ and
insures that we only need to take the orbits in $\bar{\cC^1}$.
\ere
\section{Generalizations}\label{section:generalizations}
The techniques we have developed in this paper so far can be used in the more general setting of concrete categories.
In this section we merely outline this generalization since the translations are rather straightforward and we are mainly interested in group theory. Our main reference for (concrete) categories is~\cite{AdaHerStr90}.

\mn
We start by recalling a few definitions from category theory.
For a category $\bC$, we denote the collection of objects $\obj(\bC)$ and for $A,B\in \obj(\bC)$ we let $\hom_\bC(A,B)$ denote the collection of $\bC$-morphisms from  $A$ to $B$

\bde \label{def:concrete} A  {\em concrete category} (over $\Set$) is a category $\bC$ equiped with a  functor $F:\bC\to \Set$ that is faithful, that is, for objects $A, B\in \obj(\bC)$, the hom-set restrictions
 $$F\colon \hom_\bC(A,B)\to \hom_\Set(F A, F B)$$
 are injective.
The functor $F$ is called the {\em forgetful functor}.
In~\cite{AdaHerStr90} a concrete category over $\Set$ is called a construct.
\ede

Since a given category might be concrete over several categories, the functor $F$ is part of the definition of a concrete category. 
For each $A\in \obj(\bC)$ we shall write $|A|=F(A)$.
Moreover, using that $F$ is faithful, we shall abusively identify each $f \in \hom _\bC(A,B)$ with the set function $F(f) \in \hom _{\Set}(|A|,|B|)$. Note however, that in general not all functions $f:|A|\to |B|$ correspond to morphisms of $\bC$.

\bde \label{def:embedding}
A morphism $A\stackrel{f}{\to}B$ is called an {\em embedding} if 
\begin{enumerate}
\item the underlying function $|A|\stackrel{f}{\to}|B|$ is injective, and 
\item\label{initial} $f$ is initial, i.e.~for any $C\in\obj(\bC)$, a function $|C|\stackrel{g}{\to}|A|$ underlies a $\bC$-morphism whenever  $|C|\stackrel{f\circ g}{\to}|B|$ underlies a $\bC$-morphism.
\end{enumerate}
\ede

\bre\label{rem:inverses} Note that condition~\ref{initial}.~ensures that if you can compose a morphism with $f^{-1}$ in $\Set$ then you can do so in $\bC$. More precisely, suppose that $h\in \hom_\bC(C,B)$ is so that $h(|C|)\subseteq f(|A|)$. We can then define the composition of functions $g=f^{-1}\circ h \in \hom_{\Set}(|C|,|A|)$. Now since $f\circ g=h$ and $h$ is a $\bC$ morphism, so is $g$.\ere

We now consider a simplicial complex    $X=(V,\Sigma)$ and  a concrete  category
$(\bC, F)$. 
The simplicial complex $X$ gives rise to its poset category $\bX$ in which $\obj(\bX)=\Sigma$ and 
$$\hom_\bX(A,B)=\left\{\begin{array}{ll}\{"A\le B"\}  &\mbox{ if } A\le B\\
\emptyset &\mbox{ if } A\not\le B\end{array}\right.$$
In this setting, an amalgam (cf. Definition~\ref{dfn:amalgam}) is a functor.
\bde\label{dfn:categorical simplicial amalgam}
 A {\em simplicial amalgam over $\bX$ in the (concrete) category $\bC$} is a
contravariant functor $\cG: \bX \to \bC$ such that for any pair of
simplices $\sigma, \tau$ with  $\sigma\le \tau$, the map 
$\cG("\sigma\le \tau")\colon\cG(\tau) \to \cG(\sigma)$ is an embedding (so in particular $\cG("\sigma\le \tau")$ is a monomorphism). By analogy with the group theoretic setting we shall
denote $\cG(\sigma)$ by $\cG_\sigma$ and $\cG("\sigma\le \tau")$ by $\varphi^\sigma_\tau$ and write $\cG=(\cG_\bullet,\varphi^\bullet_\bullet)$. 
\ede
\noindent Since $\bX$ and $\bC$ are fixed in this section, we shall simply call $\cG$ an amalgam.
\bde
For an amalgam $\cG$ and each $\sigma,\tau\in \Sigma$ with $\sigma\le \tau$ we define  
$$|\bar\cG_{\sigma, \tau}|= \phi^\sigma_\tau(|\cG_\tau|)\subseteq |\cG_\sigma|.$$ Note that this is a set, not an object of $\bC$.
Given an amalgam ${}^0\cG$, {\em an amalgam of type} ${}^0\cG$ is an amalgam $\cG$ such that 
 $\cG_\sigma={}^0\cG_\sigma$ for all $\sigma\in \Sigma$ and 
 $|\bar{{}^0\cG}_{\sigma, \tau}|=|\bar{\cG}_{\sigma, \tau}|$ for all 
  $\sigma,\tau\in \Sigma$ with $\sigma\le \tau$.
We call $\cG$ {\em normalized} if for any simplex $\tau$ we have
$\cG("{\bar{\tau}}\le \tau")=^0\!\cG("{\bar{\tau}}\le \tau")$.
\ede

\noindent In this setting, a morphism of amalgams (see Definition~\ref{dfn:CT hom}) is a natural transformation.
\bde We define a {\dfn homomorphism} between the amalgams
$\cG^{(1)}=\{\cG_\bullet^{(1)},  {}^{(1)}\varphi_\bullet^\bullet\}$ and
$\cG^{(2)}=\{\cG_\bullet^{(2)}, {}^{(2)}\varphi_\bullet^\bullet\}$ to be a 
 natural transformation $\phi\colon \cG^{(1)}\to\cG^{(2)}$.
That is, $\phi=\{ \phi_\sigma\mid \sigma\in\sS\}$ where $\phi_\sigma\in \hom_\bC(\cG_\sigma^{(1)},  \cG^{(2)}_{\sigma})$ is such
that, when $\sigma\le \tau$ we have 
\beq\xymatrix{
\cG^{(1)}_\sigma \ar[r]^{\phi_\sigma}& \cG^{(2)}_\sigma\\
\cG^{(1)}_\tau \ar[u]^{{}^{(1)}\varphi^{\sigma}_{\tau}} \ar[r]^{\phi_\tau} &\cG^{(2)}_\tau \ar[u]_{{}^{(2)}\varphi^{\sigma}_{\tau}}\\
}\eeq
We call $\phi$ an {\dfn isomorphism} of amalgams if it is a natural isomorphism, that is, 
if $\phi_\sigma$ is an isomorphism for all $\sigma\in \sS$. \ede

\bde\label{dfn:amal}
The simplicial amalgams over $\bX$ in $\bC$ and 
 homomorphisms of such amalgams form the objects and morphisms of a category that we shall denote 
  $\Amal_\bX(\bC)$.
It is a subcategory of the functor quasi-category $[\bX,\bC]$, which itself is (isomorphic to) a category since $\bX$ is a small category
(see~\cite{AdaHerStr90}).
\ede

\noindent One can now easily reprove Proposition~\ref{prop:amalgam normalization}.
Next, we define what is the coefficient system associated to an amalgam $\cG=\{\cG_\bullet , \psi_\bullet^\bullet\}$. 
Recall that an automorphism of $A\in \obj(\bC)$ is an isomorphism in $\hom_\bC(A,A)$. The collection of all automorphisms of $A$ is denoted $\Aut_\bC(A)$ and forms a group.

\bde\label{dfn:categorical coefficient system}
The coefficient system associated to the amalgam $\cG=\{\cG_\bullet , \psi_\bullet^\bullet\}$ is the (covariant) functor
 $\cA\colon \bX\to \Group$, where
for each simplex $\sigma\in\Sigma=\obj(\bX)$ we have  
 $$\cA(\sigma)=\{g\in \Aut_\bC(\cG_\sigma)\mid g(|\bar \cG_{\sigma, \tau}|)=|\bar \cG_{\sigma, \tau}|\mbox{ for all }\tau\mbox{ with }\sigma\le \tau\},$$
and 
for every $\sigma,\tau\in \Sigma$ with $\sigma\le\tau$ we have
$$\begin{array}{rl}
\cA("\sigma\le \tau")  \colon A_\sigma&\to A_\tau\\
f&\mapsto (\psi^\sigma_\tau)^{-1}\circ f \circ \psi^\sigma_\tau. 
\\
\end{array}
$$
As before we shall write $\cA=\{A_\bullet, \alpha_\bullet^\bullet\}$, where
 $A_\sigma=\cA(\sigma)$ and $\alpha^\sigma_\tau=\cA("\sigma\le \tau")$.
\ede
\bre
Since morphisms of amalgams are required to be embeddings, it follows from Remark~\ref{rem:inverses} that the $\alpha$'s are well-defined.
Moreover, since $((\psi^\sigma_\tau)^{-1}\circ f \circ \psi^\sigma_\tau)^{-1}=(\psi^\sigma_\tau)^{-1}\circ f^{-1} \circ \psi^\sigma_\tau$, $\alpha^\sigma_\tau$ maps $A_\sigma$ to $A_\tau$.
\ere
It is now immediate to prove Propositions~\ref{prop:amalgams to cycles} and \ref{prop:isomorphism vs cohomologous} in this categorical setting and  Theorem~\ref{mainthm2} follows.

Of course one can apply these generalizations to many natural categories such as {\bf Ring}, {\bf Top} and so on. We choose to exemplify with another category of group theoretic flavour. We refer the reader to \cite{RonSmi85} for details. 

\subsection{Sheaves on geometries.}
We shall define the notion of a sheaf on a geometry in the sense of Ronan and Smith (cf.~\cite{RonSmi85,RonSmi86}) and show that it is an example of an amalgam over a simplicial complex.
 
Let $k$ be a field and let $G$ be a finite group.
Let $\Gamma$ be an incidence geometry with type set $I$ on which  $G$ acts flag transitively. 
Let $X(\Gamma)=(V(\Gamma),\Sigma(\Gamma))$ be the 
 simplicial complex whose vertices are the objects of $\Gamma$ and whose simplices are the flags of $\Gamma$.
We let $G$ act on $X(\Gamma)$ in accordance with its action on $\Gamma$ itself.
\bde\label{dfn:sheaf}
A {\em ($G$-equivariant) sheaf} on $X(\Gamma)$ is a collection 
$$\cF_\Gamma=\{\cF_\sigma,\phi^\sigma_\tau\mid \sigma\sbe\tau\in \Sigma(\Gamma)\},$$
where $\cF_\sigma$ is a $k$-vector space and 
$\phi^\sigma_\tau\colon \cF_\tau \to \cF_\sigma$ is a linear map whenever $\sigma\sbe\tau\in\sS$. Moreover,  $\phi^\sigma_\rho\after\phi^\rho_\tau=\phi^\sigma_\tau$ whenever $\sigma\le \rho\le \tau$. Finally we require that $G$ acts on the set $\prod _\sigma \cF_\sigma$.  More precisely,  $g\in G$ acts  linearly via $\tilde g_\sigma: \cF_\sigma \to \cF_{\sigma g}$ and the following diagram commutes.
\begin{eqnarray}\label{eqn:sheaf}
\xymatrix{
\cF_\sigma \ar[r]^{\tilde g_\sigma}  & \cF_{\sigma g} \\
\cF_\tau\ar[u]^{\phi^\sigma_\tau}\ar[r]^{\tilde g_\tau}&\cF_{\tau g}\ar[u]_{\phi^{\sigma g}_{\tau g} } }. 
\end{eqnarray}
\ede

Suppose moreover that $G$ is transitive on the chambers of the geometry $\Gamma$.
In this case, much of the information contained in a sheaf is redundant. For more details on the construction see Section 2 of \cite{RonSmi86}.
Indeed, let us take $c$ a chamber of the geometry and for each subflag $\sigma$ of $c$ consider $G_\sigma$, the stabilizer of $\sigma$. 
%
A {\em stalk at $c$} is a system 
$$\cF=\{\cF_\sigma,\phi^\sigma_\tau\mid \sigma\sbe\tau\sbe c\},$$
where $\cF_\sigma$ is a $kG_\sigma$-module and $\phi^\sigma_\tau$ is a $kG_\sigma$-module  homomorphism whenever $\sigma \sbe \tau\sbe c$.

We then have (see theorem 2.3 of \cite{RonSmi86})
\bth Every stalk is the restriction (to the faces of $c$) of a unique sheaf.\eth

We shall now consider the category $\GRep$ of pairs $(G,V)$ where $G$ is a (finite) group and $V$ is a $k$ representation of $G$.  A homomorphism in $\hom_{\GRep}((G_1,V_1),(G_2,V_2))$ is a  pair $(\phi,\psi)$, where $\phi:G_1\to G_2$ is a group homomorphism, $\psi:V_1\to V_2$ is a linear map, and for all $g\in G_1$ and $v\in V_1$,  we have
 $$\psi(gv)=\phi(g)\psi(v).$$  

The functor $F\colon\GRep\to \Set , F((G,V))=G\times V$ and $F((\phi,\psi))=\phi\times \psi$ makes $(\GRep,F)$ into a concrete category. 
Let $X$ be the simplex of faces of $c$.
We then can define $\Amal_X(\GRep)$ to be the category of   simplicial amalgams over $X$ in $\GRep$. There is an obvious forgetful functor $\cG: \Amal_X(\GRep)\to \Amal_X(\Group)$ from this category into the category of  simplicial amalgams over $X$ in the category $\Group$. The following is now quite clear.
\bth
The category of sheaves over $\Gamma$ is equivalent to the fiber category $\cG_\cA$ where $\cA$ is the amalgam of parabolics of $G$ defined by $\Gamma$.
\eth

This suggests that the category $\Amal_X(\GRep)$ would be interesting to study as the various fiber categories will say something about the representation theory of various groups.


\begin{thebibliography}{1}

\bibitem{AdaHerStr90} 
J.~Adamek, H.~Herrlich, and G.~Strecker.
\newblock Abstract and concrete categories.
\newblock Aug 1990.

\bibitem{Bas93}
H.~Bass.
\newblock Covering theory for graphs of groups.
\newblock {\em Journal of Pure and Applied Algebra}, 89(1-2):3--47, October
  1993.


\bibitem{BenGraHofShp03}
C.~D.~Bennett, R.~Gramlich, C.~Hoffman, and S.~Shpectorov.
\newblock Curtis-{P}han-{T}its theory.
\newblock In {\em Groups, combinatorics \& geometry (Durham, 2001)}, pages
  13--29. World Sci. Publishing, River Edge, NJ, 2003.

\bibitem{BeSh2004}
C.~D. Bennett and S.~Shpectorov.
\newblock A new proof of a theorem of {P}han.
\newblock {\em J. Group Theory}, 7(3):287--310, 2004.

\bibitem{BloHof09a}
R.~Blok and C.~Hoffman.
\newblock A classification of Curtis-Tits amalgams.
\newblock arXiv:0907.1388v2

\bibitem{BloHof09b}
R.~J. Blok and C.~Hoffman.
\newblock Curtis-Tits groups generalizing Kac-Moody groups of type $\widetilde{A}_n$.
\newblock arXiv:0911.0824v2.

\bibitem{BloHof2011}
R.~J. Blok and C.~Hoffman.
\newblock Bass-{S}erre theory and counting rank two amalgams.
\newblock {\em J. Group Theory}, 14(3):389--400, 2011.

\bibitem{Che1995}
A.~Chermak.
\newblock Triangles of groups.
\newblock {\em Trans. Amer. Math. Soc.}, 347(11):4533--4558, 1995.

\bibitem{Cor1992}
J.~M. Corson.
\newblock Complexes of groups.
\newblock {\em Proc. London Math. Soc. (3)}, 65(1):199--224, 1992.

\bibitem{Cor1995}
J.~M. Corson.
\newblock Groups acting on complexes and complexes of groups.
\newblock In {\em Geometric group theory ({C}olumbus, {OH}, 1992)}, volume~3 of
  {\em Ohio State Univ. Math. Res. Inst. Publ.}, pages 79--97. de
  Gruyter, Berlin, 1995.

\bibitem{DelGolSte1985}
A.~Delgado, D.~Goldschmidt, and B.~Stellmacher.
\newblock {\em Groups and graphs: new results and methods}, volume~6 of {\em
  DMV Seminar}.
\newblock Birkh\"auser Verlag, Basel, 1985.
\newblock With a preface by the authors and Bernd Fischer.


\bibitem{Dun2005}
J.~Dunlap.
\newblock Uniqueness of Curtis-Phan-Tits amalgams.
\newblock Ph.D. thesis, Bowling Green State University, 2005.

 \bibitem{Gol80}
D.~M. Goldschmidt.
\newblock Automorphisms of trivalent graphs.
\newblock {\em The Annals of Mathematics}, 111(2):377--406, 1980.

\bibitem{Gr2004}
R.~Gramlich.
\newblock Weak {P}han systems of type {$C\sb n$}.
\newblock {\em J. Algebra}, 280(1):1--19, 2004.


\bibitem{GraHofShp03}
R.~Gramlich, C.~Hoffman, and S.~Shpectorov.
\newblock A {P}han-type theorem for {${\rm Sp}(2n,q)$}.
\newblock {\em J. Algebra}, 264(2):358--384, 2003.

\bibitem{Gra}
R.~Gramlich.
\newblock Developments in finite {P}han theory.
\newblock {\em Innov. Incidence Geom.}, 9:123--175, 2009.


\bibitem{Hae1991}
A.~Haefliger.
\newblock Complexes of groups and orbihedra.
\newblock In {\em Group theory from a geometrical viewpoint ({T}rieste, 1990)},
  pages 504--540. World Sci. Publ., River Edge, NJ, 1991.

\bibitem{Hae1992}
A.~Haefliger.
\newblock Extension of complexes of groups.
\newblock {\em Ann. Inst. Fourier (Grenoble)}, 42(1-2):275--311, 1992.


\bibitem{Kee2009}
P.~Keen.
\newblock In consideration of a geometry.
\newblock Master's thesis, University of Birmingham, 2009.

\bibitem{RonSmi85}
M.~A. Ronan and S.~D. Smith.
\newblock Sheaves on buildings and modular representations of Chevalley groups.
\newblock {\em J. Algebra}, 96:319--346, 1985.

\bibitem{RonSmi86} 
M.~A. Ronan and S.~D. Smith.
\newblock Universal presheaves on group geometries, and modular representations.
\newblock {\em J. Algebra}, 102:135--154, 1986.

\bibitem{Se1980}
J.-P. Serre.
\newblock {\em Trees}.
\newblock Springer-Verlag, Berlin, 1980.

\bibitem{Sta1991}
J.~R. Stallings.
\newblock Non-positively curved triangles of groups.
\newblock In {\em Group theory from a geometrical viewpoint ({T}rieste, 1990)},
  pages 491--503. World Sci. Publ., River Edge, NJ, 1991.

\bibitem{Wie2007}
C.~Wiedorn.
\newblock Flag-transitive {$c$}-extensions of the {$F_4(2)$}-building---the
  amalgams.
\newblock {\em European J. Combin.}, 28(2):641--652, 2007.


\end{thebibliography}
\end{document}